# Dynamics of a hybrid vibro-impact oscillator: canonical formalism

Maor Farid

*Massachusetts Institute of Technology, 77 Massachusetts Ave., Cambridge, MA 02139, United States*
*Faculty of Mechanical Engineering, Technion – Israel Institute of Technology, Haifa 3200003, Israel*
*faridm@mit.edu*

**Abstract**

Hybrid vibro-impact (HVI) oscillations is a strongly nonlinear dynamical regime that involves both linear oscillations and collisions under periodic, impulsive, or stochastic excitation. This regime arises in various engineering systems, such as mechanical components under tight rigid constraints, seismic-induced sloshing in partially-filled liquid storage tanks, and more. The adaptive nonlinearity of the HVI oscillator is used by the HVI-nonlinear energy sink as an effective vibration mitigation solution for broad energy and frequency range. Due to the extreme nonlinearity of this regime, traditional analytical methods are inapplicable for the description of its transient dynamics. In the current work, we model the HVI oscillator by a forced particle in a truncated quadratic potential well with infinite depth. The slow flow dynamics of the system in the vicinity of primary resonance is described by canonical transformation to action-angle (AA) variables and the corresponding reduced resonance manifold (RM). Two types of bifurcation are examined. The former is associated with transition between linear oscillations and the HVI-regime and vice versa, and the latter with reaching a chosen maximal transient energy level. The transition boundaries on the forcing parameters plane associated with both bifurcation types are obtained analytically. The maximal transient energy level obtained for any given set of forcing parameters is described analytically as well. The energy jumps associated with the bifurcation of type I and crossing the corresponding transition boundary are obtained. Two underlying dynamical mechanisms that govern the occurrence of bifurcations are identified. They correspond to two distinct scenarios: in the first scenario, the energy of the slow flow gradually reaches the threshold energy level and is thus referred to as the "maximum" mechanism. The second, potentially more dangerous scenario, involves abrupt transitions of the system's energy response from a relatively small value to the threshold energy level. This pattern is related to the passage of the slow-flow phase trajectory through the saddle point of the RM, and thus is referred to as the "saddle" mechanism. Both mechanisms are universal for systems that undergo escape from a potential well. All theoretical results are in complete agreement with full-scale numerical simulations.

## 1. Introduction

When various engineering systems are exposed to external disturbances, they might develop strongly nonlinear behavior that involves collisions, impacts, and chattering. Examples include machine and structural components subjected to kinematic constraints such as nuclear fuel rods subjected to flow-induced vibration [1, 2], strongly-nonlinear sloshing regimes in partially-filled liquid storage tanks under seismic excitation [3, 4, 5, 6, 7, 8], and more. Those systems can be linear in their nature, but the existence of impacts serves as a source of non-smoothness and the strongest non-linearity possible. Depending on the magnitude of the excitation, the system may develop multiple dynamical regimes, such as small-amplitude linear oscillations, and sustained collisions separated by linear motion between each consequent collisions. Hence, the latter is referred to as the hybrid vibro-impact (HVI) regime. In most cases, this intensive regime can lead to substantial accelerations, stresses, and finally increased wear of the system, which might lead to catastrophic consequences. Moreover, analysis and prediction of the system's response and resistance become substantially more difficult when in most cases the underlying equations are unsolvable. On the other hand, the HVI-regime is used for engineering purposes such as vibration mitigation. The HVI nonlinear energy sink (NES) [9] is a nonlinear passive energy absorber (PEA) that takes advantage of the coexistence of both linear and strongly nonlinear regimes to yield an adaptive nonlinearity that changes its dynamical properties in accordance with the magnitude and frequency of the external excitation. For low-energy excitation, the linear regime takes place, and the HVI-NES is practically equivalent to the well-known tuned mass damper (TMD) [10, 11]. For high-energy excitation, the nonlinearity of the HVI-NES is switched-on and HVI-regime takes place. In this case, the hybrid absorber acts as a generalization of both the TMD and the vibro-impact (VI) NES [12, 13]. The NESs are known for their vibration mitigation capabilities for wide frequency ranges. They perform ideally when a transient resonance capture (TRC) [14] with the primary structure takes place. However, when the primary structure is exposed to low-energy excitation and TRC is not achieved, and therefore the NES performs poorly. Contrarily, the TMD performs well only when sustained resonance takes place and the excitation is of relatively low magnitude, otherwise, it loses its linear properties. The HVI-NES is a hybridization or generalization of both TMD and VI-NES and thus exhibits enhanced absorption characteristics for broader energy and frequency ranges in comparison to both the TMD and well-known NES designs. Often, the energy absorption capabilities of PEAs are evaluated by the maximal instantaneous energy measure captured by the absorber. Due to the transient nature of the TRC mechanism, this quantity cannot be evaluated based on the steady-state response of the system and thus, it heavily relies on numerical simulations. However, the existence of adaptive and non-smooth nonlinearity poses a mathematical challenge in describing and predicting the dynamical responses which are governed by transient phenomena under external disturbances. Hence, traditional asymptotic methods such as perturbation-based approaches (multiple scales, Lindstedt–Poincaré [15], complexification averaging [16, 17], harmonic balance [18]) become inapplicable. Moreover, those methods describe the response of forced systems in steady-state, mainly in conditions of primary resonance, and allow to obtain safety boundaries in the excitation space, but don't provide insights and understanding regarding the transient dynamics of the system. In other terms, the existence of vibro-impacts contributes a hardening nonlinearity that leads to an increase in the response frequency. Therefore, most of the methods mentioned above are hardly applicable. In this paper, we focus on the analytical description of the transient dynamics of an HVI-oscillator under external harmonic forcing. More specifically, we model the hybrid oscillator as a classical particle in one dimensional truncated quadratic potential well with infinite depth. The truncation of the equivalent potential well allows the coexistence of both linear oscillations and the HVI-regime. The oscillator is initially at rest at the bottom of the well when the external force is switched on. It was shown in previous works [19, 20, 21, 22] that transient dynamical regimes can be accurately described using the reduced resonance manifold (RM) of the system that describes the slow dynamics of the system. Distinct dynamical regimes and their corresponding bifurcations are related to the topological modifications of a special phase trajectory on the RM that correspond to a set of initial conditions. The trajectory that corresponds to zero initial conditions is referred to as the limiting phase trajectory (LPT) [23, 24].

In the current work, we use canonical transformation to the action-angle (AA) variables and a proper averaging technique to describe all dynamical regimes that arise in the HVI-oscillator system and formulate the corresponding transition boundaries in the forcing parameters space. The main goal of the current work is to provide an analytical prediction of the particle's dynamical response and its maximal transient energy level (or absorption rate when thought as a vibration mitigation solution) under monochromatic harmonic forcing. Moreover, we aim to describe all possible dynamical regimes that arise in the system, the boundaries associates with the transitions from one regime to another, and the underlying dynamical mechanism that govern those transitions.

This paper is structured as follows: In Section 2 the dynamical model of the HVI-oscillator is described. In Section 3 AA formalism is performed and the dynamical regimes are explained and described analytically. The resonance manifold is computed and its structure is investigated. In Section 4 we describe the underlying bifurcation mechanisms that govern the transitions between dynamical regimes and energy levels, and obtain the corresponding transition boundaries on the forcing parameters space. Frequency response curves are obtained, as well as a mapping between forcing parameters and the resulting maximal transient energy. Section 5 includes numerical validations of the analytic results. Section 6 is devoted to the concluding remarks.

## 2. Model description and action-angle formalism

The system considered is a HVI-oscillator with mass $m$ moving in a channel of length $2d$, with elastic collisions at the ends of the channel, as shown in Fig. 1. The oscillator is subjected to a time-dependant forcing $\bar{F}(t)$ and is coupled to the rigid walls by a pair of linear springs of stiffness $k/2$ each. The displacement of the oscillator with respect to its equilibrium position is denoted by $x$.

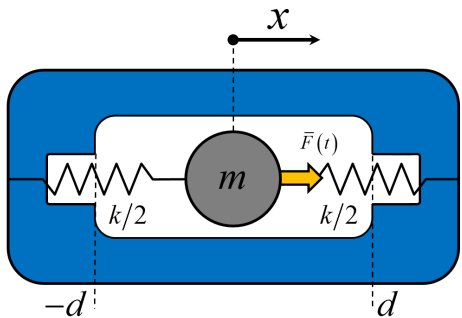

Figure 1: Scheme of the forced hybrid vibro-impact oscillator

The normalized equation of motion of the system in shown in Eq. (1), where $q = x/d$ and $\tau = \sqrt{k/m}t$ are normalized displacement and time scale, respectively.

$$\ddot{q} + q + \kappa \sum_j \dot{q}(\tau_j^-)\delta(\tau - \tau_j) = F\cos(\Omega\tau) \qquad (1)$$

Here $F$ and $\Omega$ are the non-dimensional forcing amplitude and frequency, respectively. Dot denote differentiation with respect to the non-dimensional time-scale $\tau$, $\tau_j$ is the time instance of the $j^{th}$ collision since the forcing switched-on at $\tau = 0$, parameter $\kappa$ is the restitution coefficient, and $\delta(\tau)$ is the Dirac delta function. We consider purely-elastic collisions and take a restitution coefficient of unity $\kappa = 1$ accordingly. Thus, since the system is conservative, the motion of the HVI-oscillator can be modeled by a truncated quadratic potential well with infinite depth as shown in Eq. (2) and Fig. 2. The transition between linear oscillations and vibro-impact regime corresponds to $|q| = 1$, i.e. instantaneous energy of $E = 1/2$.

$$U(q) = \begin{cases} \frac{1}{2}q^2 & , 0 < E \leq 1/2 \\ \infty & , E > 1/2 \end{cases} \tag{2}$$

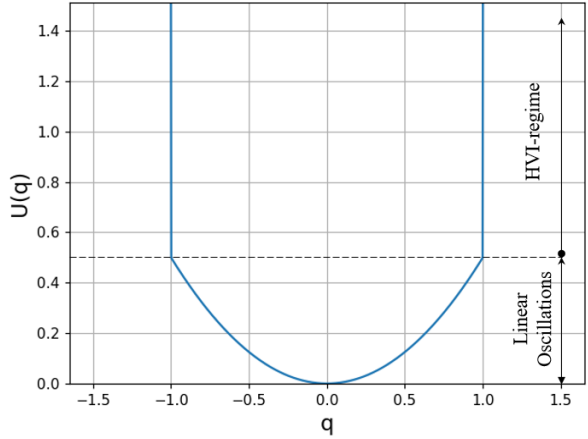

Figure 2: The equivalent potential well of a hybrid vibro-impact oscillator. Energy values of $U(q) \in (0, 1/2)$ correspond to the regime of linear oscillations and values of $U(q) \in [1/2, \infty)$ correspond to the HVI-regime. The limiting case of $U(q) = 1/2$ corresponds to a bifurcation of type I (dashed line).

## 3. Dynamical regimes and resonance manifold topology

In the current section, the dynamical regimes of the HVI-oscillator are described and analyzed using analytical tools. Moreover, the transition boundaries in the space of forcing parameters are formulated and the underlying mechanisms that govern the transitions are explored.

*3.1. Linear oscillations*

Although the system is linear and hence solvable for $E \in (0, 1/2)$, we are not in its steady-state solution, but specifically in the conditions that lead to the transient rise in the response energy of the forced HVI-oscillator for modification of the forcing parameters. For the sake of approximate analysis in the vicinity of primary resonance, the following small-amplitude forcing and detuning parameters are adopted:

$$F = \epsilon f, \quad \Omega = 1 + \epsilon \sigma \tag{3}$$

Here $\epsilon \ll 1$ is a small parameter. Then, following [22], we obtain the approximate relation between detuning parameter $\sigma$, the maximal averaged energy level of the forced linear oscillator $\xi$, and the resulting critical forcing amplitude $f$:

$$f_m(\sigma, \xi) \approx \sqrt{2\xi}|\sigma| \tag{4}$$

The inverse relation of Eq. (4), i.e. $\xi(\sigma, f_m)$, defines iso-energy contours in the forcing parameters plane, as shown in Fig. 3. The condition for transition between linear oscillations and vibro-impact regime corresponds to the transition boundary in the forcing parameters space, which is associated with the limiting dynamical regime of linear oscillations of amplitude of unity and maximal instantaneous energy of $\xi = 1/2$.

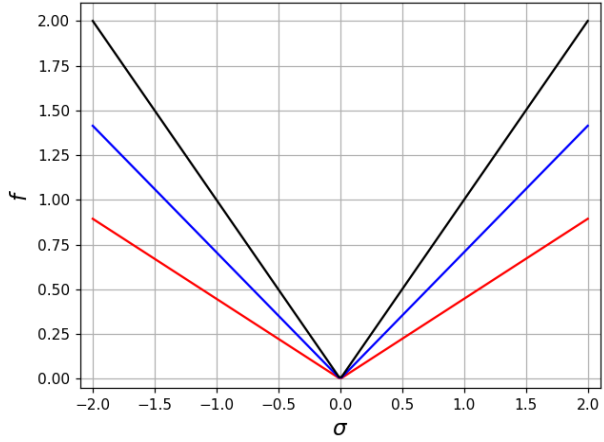

Figure 3: Critical forcing amplitude versus the detuning parameter $f_m(\sigma|\xi)$ or transition boundaries that correspond to bifurcations of type II, for given values of maximal averaged transient energy: $\xi = 0.1, 0.5, 1$ in red, blue, black, respectively. The latter corresponds to bifurcation of type I.

The minima in Fig. 3 correspond to the natural frequency of the HVI-oscillator. The sharp shape of the curve is universal for systems that undergo escape from a potential well [25, 26, 27] and arises in various fields, such as chemistry, physics and engineering [28, 29, 30, 31, 32, 33].

3.2. Hybrid vibro-impact regime

When the HVI-regime takes place, the equation of motion of the system in Eq. (1) demonstrates the strongest possible non-linearity and non-smoothness. Hence, perturbation-based approaches are inapplicable in the current case. The generalized basis functions that describe the motion of the HVI oscillator $g(\tau), g'(\tau)$ are given in Eq. (5) and shown in Fig. 4. Functions $\bar{\tau}(\tau)$ and $\bar{e}(\tau)$ are the basis functions that correspond to the VI-oscillator [34, 35], and the basis functions that correspond to the regime of linear oscillations are the standard trigonometric functions $\sin \tau$ and $\cos \tau$. The basis functions introduced serve as a generalization of both sets of basis functions and controlled by parameter $\beta$ which is determined by the energy of the oscillator according to Eq.(5). The limiting values $\beta = \pi/2$ ($E = \pi^2/8$) and $\beta = 0$ ($E = 1/2$) correspond to pure linear oscillations and the VI motion, respectively, and $\beta \in (0, \pi/2)$ corresponds to the HVI regime. The Fourier expansion of the generalized basis function $g(\tau)$ is shown in Eq. (6).

$$g(\tau) = \frac{\sin(\beta \bar{\tau}(\tau))}{\sin \beta}, \quad g'(\tau) = \beta \frac{\cos(\beta \bar{\tau}(\tau))}{\sin \beta} \bar{e}(\tau), \quad E = \frac{1}{2}\left(\frac{\beta}{\sin \beta}\right)^2$$
$$\bar{\tau}(\tau) = \frac{2}{\pi}\arcsin(\sin \tau), \quad \bar{e}(\tau) = \frac{d\bar{\tau}(\tau)}{d\tau} = \text{sgn}(\cos \tau) \tag{5}$$

$$g(\tau) = \sum_{n=1}^{\infty} b_n \sin(n\tau), \quad b_n = \frac{4\beta}{\pi(\beta^2 n^2 - 1)}\cot\left(\frac{\pi}{2\beta}\right)\sin\left(\frac{\pi}{2}n\right) \tag{6}$$

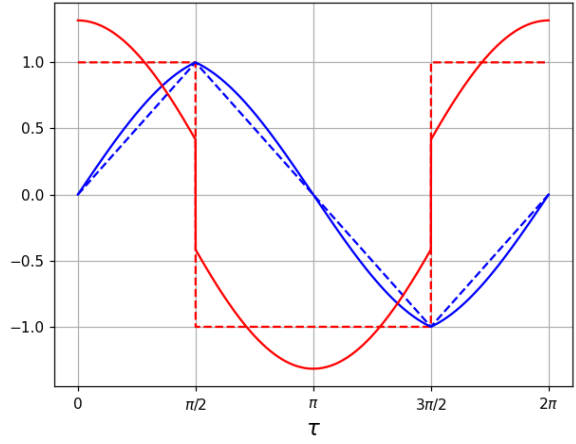

Figure 4: Generalized basis functions for the HVI-oscillator, $g(\tau), g'(\tau)$ for $\beta = 1.25$ in solid blue and red lines, respectively, and the basis functions of the VI-oscillator $\bar{\tau}(\tau), \bar{e}(\tau)$ (correspond to $\beta = 0$) in dashed blue and red lines, respectively.

Due to the non-smoothness of the basis functions that correspond the HVI-regime, infinite series of Fourier term are required to approximate its response under periodic excitation. Instead of using the cumbersome infinite Fourier series in Eq. (6), the generalized basis functions in Eq. (5) can be easily adopted. Following [36, 19], we derive the Hamiltonian of the system:

$$H = H_0(p,q) - Fq\sin(\Omega\tau), \quad H_0(p,q) = \frac{1}{2}\left(p^2 + q^2\right) \qquad (7)$$

Here $p = \dot{q}$ is the momentum of the oscillator and component $H_0(p,q)$ describes the free motion of the oscillator in the potential well, i.e. corresponding to the absence of an external forcing. Therefore it is referred to as the conservative component of the Hamiltonian. The transformation to action-angle (AA) variables is obtained by the following formulas [37]:

$$I(E) = \frac{1}{2\pi} \oint p(q,E)dq, \quad \theta = \frac{\partial}{\partial I}\int_0^q p(q,I)dq \qquad (8)$$

Here $I$ and $\theta$ are the action and angle variables, respectively. In Eq. (8), $H_0 = E$ defines a constant energy level. Inverting Eq. (8) yields explicit formulas for the canonical variables transformation, $p(I,\theta)$ and $q(I,\theta)$, and the conservative component of the Hamiltonian can be formulated only using the action variable: $H_0 = E(I)$. Then, the Hamiltonian of the system can be written in terms of AA variables as follows:

$$H = H_0(I) - Fq(I,\theta)\sin(\Omega\tau) \qquad (9)$$

Due to periodicity of $2\pi$ of the angle variable $\theta$, Eq. (9) can be reframed in terms of Fourier series [38]:

$$H = H_0(I) + \frac{iF}{2}\sum_{n=-\infty}^{\infty} q_n(I)\left(e^{i(n\theta+\Omega\tau)} - e^{-i(n\theta-\Omega\tau)}\right), \quad q_n = \bar{q}_{-n} \qquad (10)$$

Here bar represents the complex conjugate. Then, the Hamilton equation takes the fol-

lowing form:

$$\begin{aligned}\dot{I} &= -\frac{\partial H}{\partial \theta} = \frac{F}{2} \sum_{n=-\infty}^{\infty} n q_n(I) \left(e^{i(m\theta+\Omega\tau)} - e^{i(m\theta-\Omega\tau)}\right) \\ \dot{\theta} &= \frac{\partial H}{\partial I} = \frac{\partial H_0}{\partial I} + \frac{iF}{2} \sum_{n=-\infty}^{\infty} \frac{\partial q_n(I)}{\partial I} \left(e^{i(m\theta+\Omega\tau)} - e^{i(m\theta-\Omega\tau)}\right)\end{aligned} \quad (11)$$

In the current work, the vicinity to primary resonance is considered. Hence, we assume slow evolution of the phase variable $\nu = \theta - \Omega\tau$ and fast variation of all other harmonics $n > 1$. Averaging over the fast phase variables yields the following slow evolution equations:

$$\begin{aligned}\dot{J} &= -\frac{F}{2}\left(q_1(J)e^{i\nu} + \bar{q}_1(J)e^{-i\nu}\right) \\ \dot{\nu} &= \frac{\partial H_0}{\partial I} + \frac{iF}{2}\sum_{n=-\infty}^{\infty}\frac{\partial q_n(I)}{\partial I}\left(e^{i(m\theta+\Omega\tau)} - e^{i(m\theta-\Omega\tau)}\right)\end{aligned} \quad (12)$$

Here $J(\tau) = \langle I(\tau) \rangle$ is the averaged action variable. Using Eq. (12), one can obtain the following conservation rule:

$$C(J,\nu) = H_0(J) - \frac{iF}{2}\left(q_1(J)e^{i\nu} - \bar{q}_1(J)e^{-i\nu}\right) - \Omega J \quad (13)$$

Eq. (13) defines a family of resonance manifolds (RMs), each corresponds to a distinct set of initial conditions. In the current study we focus on a forced HVI-oscillator that begins its motion from a rest state. The contour on the RM that correspond to this case, i.e. $H_0(J) = 0$, is referred to as the limiting phase trajectory (LPT) [17, 16, 23, 24]. The averaged action is given by the following expression:

$$J(\xi) = \frac{1}{\pi}\left(\varphi(\xi) + 2\xi \arctan\left(\frac{1}{\varphi(\xi)}\right)\right), \quad \varphi(\xi) = \sqrt{2 \cdot \max\left(\xi, \frac{1}{2}\right) - 1} \quad (14)$$

Here, as mentioned in the previous section, $\xi(\tau) = \langle E(\tau) \rangle$ is the averaged energy of the system, where averaging is performed over the fast phase variables. The relation between the oscillator energy and frequency are shown in Eq. (15) and Fig. 5.

$$\omega(\xi) = \left(\frac{\partial J}{\partial \xi}\right)^{-1} = \frac{\pi}{2}\left(\arctan\left(\frac{1}{\varphi(\xi)}\right)\right)^{-1} \quad (15)$$

The expressions in Eq. (14)-(15) are valid for both linear oscillations and HVI-regime. Detailed derivations of these expressions are presented in Appendix A. As one can see, the response frequency increased with energy, as expected in the case of hardening non-linearity associated with elastic impacts.

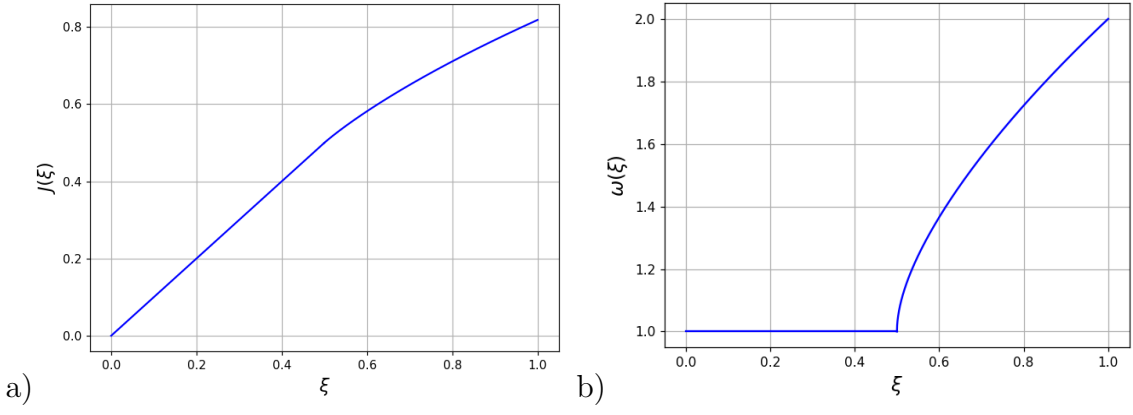

Figure 5: a) The averaged action, and b) the frequency of the HVI-oscillator vs. the averaged energy $\xi$, according to Eq. (14) and Eq. (15), respectively. Domains $\xi \in (0, 1/2)$ and $\xi \geq 1/2$ correspond to linear oscillations and HVI-regime, respectively. Monotonous increase of the frequency $\omega(\xi)$ stems from hardening non-linearity associated with elastic impacts.

In accordance to Eq. (8), the response of the oscillator is given by the following expression in terms of AA variables:

$$q(\theta, \xi) = \sqrt{2\xi} \sin\left(\frac{\theta}{\omega(\xi)}\right) \quad (16)$$

As can be easily seen, Eq. (14) is not invertible, i.e. one cannot explicitly express the averaged energy in terms of the averaged action $\xi(J)$. Consequently, all values following from the AA transformation will be parameterized through the averaged energy $\xi(t)$ instead of the averaged action $J(t)$. Specifically, the conservation law is recast in the following form:

$$C(\nu, \xi) = \xi - \frac{\epsilon f}{2} a_1(\xi) \cos(\nu) - (1+\epsilon\sigma)J, \quad a_1(\xi) = \begin{cases} \sqrt{2\xi} & \xi \in (0, 1/2) \\ \frac{\sqrt{2\xi} 2\omega^2(\xi) \sin\left(\frac{\pi}{\omega(\xi)}\right)}{\pi(\omega^2(\xi)-1)} & \xi \in [1/2, \infty) \end{cases} \quad (17)$$

Here $a_1(\xi)$ is the coefficient of the first term in the Fourier expansion of $q(J, \theta)$ in Eq. (16). The value of the first Fourier coefficient $a_1(\xi)$ is presented in Fig. 6. For high energy levels, the following limit is obtained: $\lim_{\xi \to \infty} a_1(\xi) = 4/\pi$. Detailed derivation is shown in Appendix A. As mentioned above, further analysis focuses on the level line of the phase cylinder $\nu \in [0, 2\pi]$, $\xi \in [0, \infty]$ that correspondents to zero initial conditions, i.e. the LPT. The latter is defined by the following relation $C(\nu, \xi) = 0$. The maximal averaged energy level reached by the HVI-oscillator $\xi(\sigma, f)$ corresponds to the maximal height on the phase cylinder to which the LPT reaches from the bottom of the cylinder $\xi = 0$.

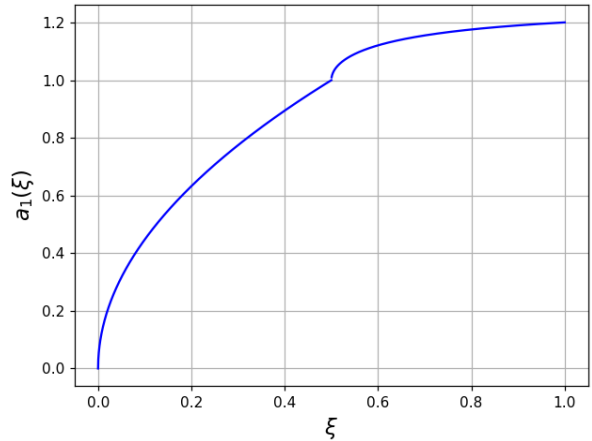

Figure 6: The first Fourier coefficient of $q(\theta, \xi)$ versus the averaged energy of the HVI-oscillator, $a_1(\xi)$.

### 3.3. Structure of the RM

Now, let us explore the topological structure (phase portrait) of the RM, defined by the conservation law in Eq. (17). We assume that the forcing is switched at $\tau = 0$ and zero-initial conditions are considered: $q(0) = \dot{q}(0) = 0$. Moreover, we assume that the oscillator is captured to the RM throughout its motion. The LPT corresponds to zero action $J = 0$, which satisfies the following expression:

$$C(\nu, \xi) = 0 \qquad (18)$$

The fixed points of the phase portrait correspond to the solutions of the following equations:

$$\frac{\partial C(\nu, \xi)}{\partial \nu} = 0, \quad \frac{\partial C(\nu, \xi)}{\partial \xi} = 0 \qquad (19)$$

The first equation yields that the stationary points can be only on lines $\nu_0 = 0, \pi$ on the phase cylinder. The solution of the second equation in Eq. (19) is awkward and can hardly be solved analytically. Thus, we explore its solutions using both numerical and approximate tools, as shown in the next section.

## 4. Bifurcations and underlying mechanisms

In the current section, the mechanisms that govern the bifurcations or transitions between distinct dynamical regimes are explored using both analytical and numerical tools. Here we define two types of bifurcations or transitions. Bifurcation of type I is associated with transition between dynamical regimes, i.e. linear oscillations and HVI-regime. Bifurcation of type II corresponds to reaching a chosen transient energy level of $\tilde{\xi}$, and therefore associated with the transition between the state of $\xi < \tilde{\xi}$ to state of $\xi \geq \tilde{\xi}$. In other terms, the bifurcation of type I is a particular case of the bifurcation of type II that corresponds to $\tilde{\xi} = 1/2$. Each of the two bifurcation types is represented by a transition boundary in the forcing parameters space $f - \sigma$. The main goal of this section is to express the transition boundaries of both bifurcation types analytically, and to predict the resulting maximal response energy obtained for any given set of forcing parameters.

In various previous studies [28, 29] it was shown that the transition boundary curves (there referred to as escape curves), share a common property of a sharp minimum. Recent works [21, 22, 39, 20] have shown that the sharp 'dip' of the transition boundary corresponds to intersection between two curves associated with two distinct underlying dynamical mecha-nisms, originated from the nonlinear features of the equivalent potential well. The former, corresponds to transition of the LPT trough a saddle point of the RM before reaching the critical energy level $\tilde{\xi}$, and hence called 'saddle mechanism'. In this mechanism, two branched of the LPT intersect at the saddle point, leading to a sudden increase or 'jump' in the response energy up to level $\tilde{\xi}$. The abrupt nature of the dynamical responses associated with this mechanism makes it potentially hazardous for engineering systems, and attractive for vibration mitigation purposes. The second mechanism, corresponds to direct motion of the LPT towards the chosen maximal energy level $\tilde{\xi}$, and hence referred to as 'maximum mechanism'. This regime is characterized by a gradual increase in the averages energy of the oscillator, and thus less dangerous for engineering systems.

*4.1. Maximum mechanism*

First, we begin with the regime of linear oscillations. In this case, the critical forcing amplitude that corresponds to a type II bifurcation for energy level of $\tilde{\xi}$ is obtained from Eq. (4) and written as follows:

$$f_m(\sigma|\tilde{\xi}) = \sqrt{2\tilde{\xi}}|\sigma|, \quad \tilde{\xi} \in \left(0, \frac{1}{2}\right) \tag{20}$$

For the HVI-regime, following Eq. (18)-(19), the maximum mechanism corresponds to the following expression:

$$C(\nu = 0, \tilde{\xi}|f_m) = \tilde{\xi} - \frac{\epsilon f_m}{2} a_1(\tilde{\xi}) - (1 + \epsilon\sigma)J(\tilde{\xi}) = 0 \tag{21}$$

Using Eq. (17) and Eq. (20) the following expression is obtained:

$$f_m(\sigma|\tilde{\xi}) = \begin{cases} \sqrt{2\tilde{\xi}}|\sigma| & \tilde{\xi} \in \left[0, \frac{1}{2}\right) \\ \frac{-2}{\epsilon a_1(\tilde{\xi})} \left((1+\epsilon\sigma)J(\tilde{\xi}) - \tilde{\xi}\right) & \tilde{\xi} \in \left[\frac{1}{2}, \infty\right) \end{cases} \tag{22}$$

for the linear regime, *only* maximum mechanism can take place. This, as will be shown in the next subsection, is due to fact that the saddle point is identically equal to one half $\xi_s = 1/2$, which means that both the maximum and saddle points overlap. Bifurcations of type I ($\tilde{\xi} = 1/2$) and type II for $\tilde{\xi} = 1$ (chosen arbitrarily) are shown in Fig. 7 and Fig. 8, respectively. In Fig. 7, even though saddle mechanisms take place in a global perspective of the RM, we treat it as a maximum mechanism because the bifurcation threshold considered for type II overlooks energy values above the transition boundary, i.e. $\xi > \tilde{\xi} = 1/2$.

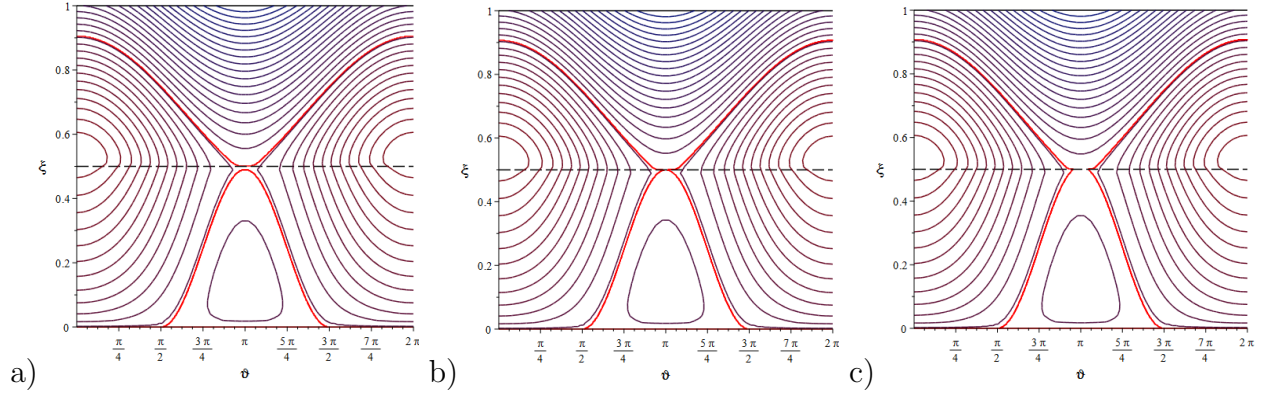

Figure 7: Bifurcation of type I through the maximum mechanism in phase portrait, defined by $C(\nu, \xi)$. The LPT is marked by a red line. Corresponds to $\tilde{\xi} = 1/2$ by definition (dashed line). For detuning of $\sigma = 1$, and a) $f = 0.99$, i.e. below the transition boundary b) on the transition boundary $f = f_m(\sigma = 1, \tilde{\xi} = 1/2) = 1.0$, c) above the transition boundary $f = 1.01$.

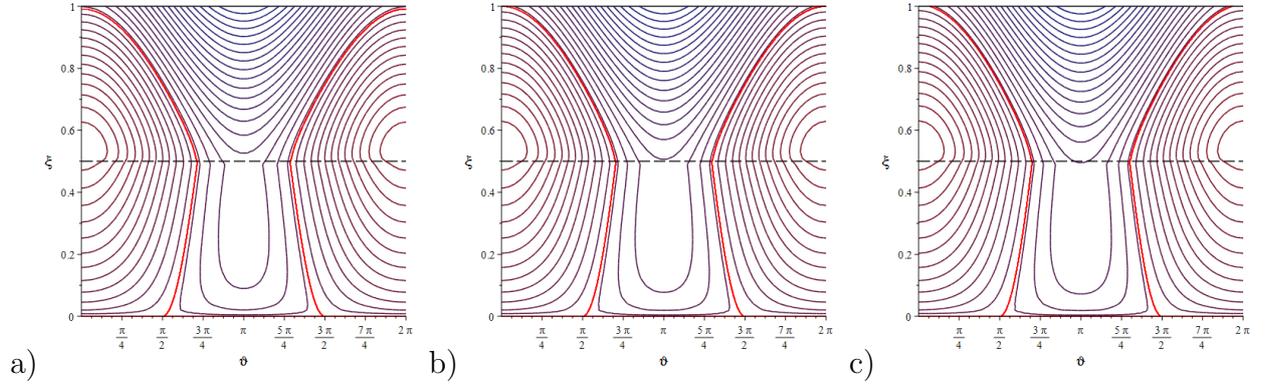

Figure 8: Bifurcation of type II through the maximum mechanism for threshold energy of $\tilde{\xi} = 1$ in phase portrait, defined by $C(\nu, \xi)$. The LPT is marked by a red line and the energy level associated with the saddle point is marked by a dashed line. For detuning of $\sigma = 1$, and a) $f = 1.6$, below the transition boundary, b) $f = f_m(\sigma = 1, \tilde{\xi} = 1) = 1.6637$, on the transition boundary, c) $f = 1.7$ above the transition boundary.

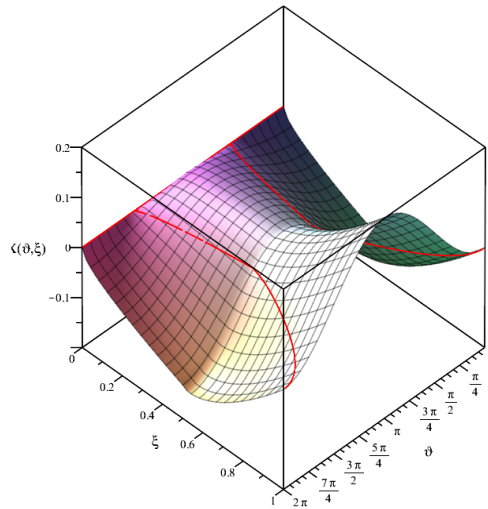

Figure 9: A 3D plot of $C(\nu, \xi)$. The LPT is marked by a red line. Bifurcation of type II through the maximum mechanism for threshold energy of $\tilde{\xi} = 1$, detuning of $\sigma = 1$, and forcing amplitude of $f = 1.6637$.

## 4.2. Saddle mechanism

Following Eq. (18)-(19), the saddle mechanism corresponds to the following relations:

$$
\begin{aligned}
C(\nu = \pi, \xi_s | f_s) &= \xi_s + \frac{\epsilon f_s}{2} a_1(\xi_s) - (1 + \epsilon\sigma) J(\xi_s) = 0 \\
\frac{\partial C}{\partial \xi}(\nu = \pi, \xi_s | f_s) &= 1 + \frac{\epsilon f_s}{2} a_1'(\xi_s) - (1 + \epsilon\sigma) J'(\xi_s) = 0
\end{aligned}
\quad (23)
$$

Here $\xi_s$ is the energy level associated with the saddle point of the RM, and $f_s$ is the critical forcing amplitude required for reaching the saddle point. After the LPT reaches the energy level associates with the saddle point $\xi_s$, it intersects with its upper branch and then 'jumps' to the critical energy level $\tilde{\xi}$. We eliminate the forcing amplitude by applying a simple algebraic manipulation on Eq. (23), to yield the following implicit relation between the detuning parameter $\sigma$ and energy values of the stationary points on line $\nu = \pi$ on the phase cylinder:

$$
\sigma = \frac{1}{\epsilon} \left( \frac{a_1(\xi_0) - a_1'(\xi_0)\xi_0}{a_1(\xi_0)J'(\xi_0) - a_1'(\xi_0)J(\xi_0)} - 1 \right)
\quad (24)
$$

Since the relation between the energy of the stationary points of the RM $\xi_0$ and $\sigma$ in Eq. (24) is implicit, let us explore it numerically. Explicit plot of Eq. (24) is presented in Fig. 10. As one can see, a degenerate saddle point exists for any value of forcing amplitude and detuning, and equals to $\xi_s = 0.5$. The degeneration of the saddle point corresponds to the non-smoothness of Eq. (24) that stems from the piece-wise definition of Eq. (14) and Eq. (17).

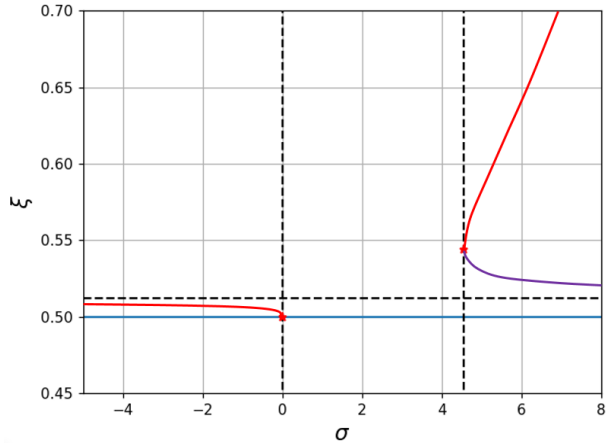

Figure 10: Plot of the stationary points of the RM in the detuning-energy plane according to the explicit expression Eq. (24), for $\epsilon = 0.1$; saddle points, minima, and maxima are colored in blue, red, and purple, respectively. Black dashed lines correspond to vertical and horizontal asymptotes, and bifurcation points $(0, 0.5), (4.5636, 0.5435)$ are marked with red stars. The horizontal asymptote corresponds to $0.51229$.

Substituting Eq. (24) into the first equation in Eq. (23) and taking $\xi_s = 1/2$, we obtains the following expression for the critical forcing amplitude associated with bifurcation of type II through the saddle mechanism:

$$
f_s(\sigma) = \frac{2}{\epsilon a_1(\xi_s)} \left( (1 + \epsilon\sigma) J(\xi_s) - \xi_s \right) = \sigma
\quad (25)
$$

From Eq. (25) we can see that the branch of the transition boundary associated with type II bifurcation overlaps the right branch of the transition boundary associates with type I bifurcation. Demonstration of the bifurcation of type II for $\tilde{\xi} = 1$ through the saddle mechanism from the perspective of the phase portrait is shown in Fig. 11.

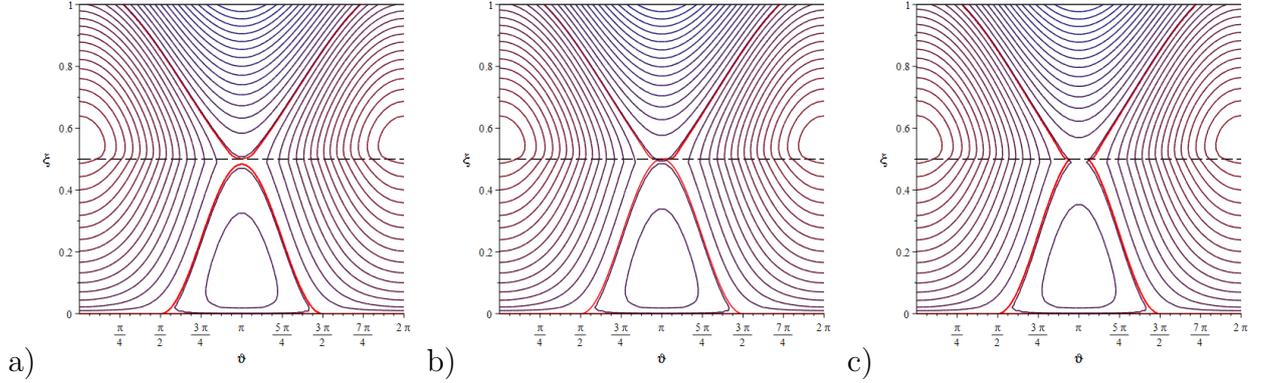

a) b) c)

Figure 11: Bifurcation of type II for $\tilde{\xi} = 1$ through the saddle mechanism in phase portrait, defined by $C(\nu, \xi)$. The LPT is marked by a red line. For detuning of $\sigma = 1.5$, and a) $f = 1.475$, i.e. below the transition boundary , b) $f = f_s(\sigma = 1.5, \tilde{\xi} = 1) = 1.5$, i.e. on the transition boundary, c) $f = 1.525$, i.e. above the transition boundary.

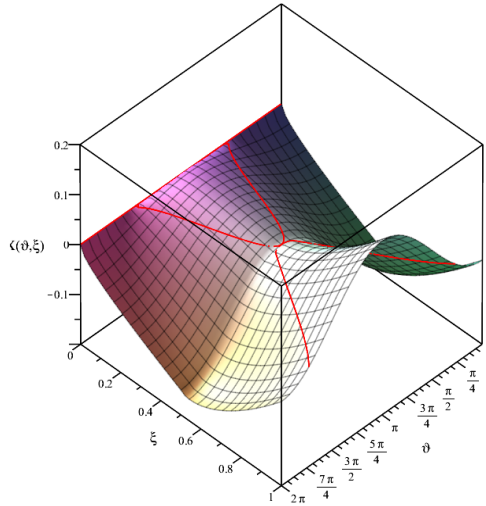

Figure 12: Bifurcation of type II through the saddle mechanism for threshold energy of $\tilde{\xi} = 1$ in a 3D plot of $C(\nu, \xi)$. The LPT is marked by a red line. For detuning and forcing amplitude of $\sigma = 1.5$ and $f = 1.5$, respectively.

Coexistence of both maximum and saddle mechanisms corresponds to the intersection of transition boundary branches, i.e. $f_m(\sigma^* | \tilde{\xi}) = f_s(\sigma^*)$, where $\sigma^*$ is the detuning value that corresponds to the coexistence of both mechanisms, and given by the following expression:

$$f^*, \sigma^*(\tilde{\xi}) = \frac{1}{\epsilon}\left(\frac{a_1(\xi_s)\tilde{\xi} + a_1(\tilde{\xi})\xi_s}{a_1(\tilde{\xi})J(\xi_s) + a_1(\xi_s)J(\tilde{\xi})} - 1\right) = \frac{1}{\epsilon}\left(\frac{a_1(\tilde{\xi}) + 2\tilde{\xi}}{a_1(\tilde{\xi}) + 2J(\tilde{\xi})} - 1\right) \qquad (26)$$

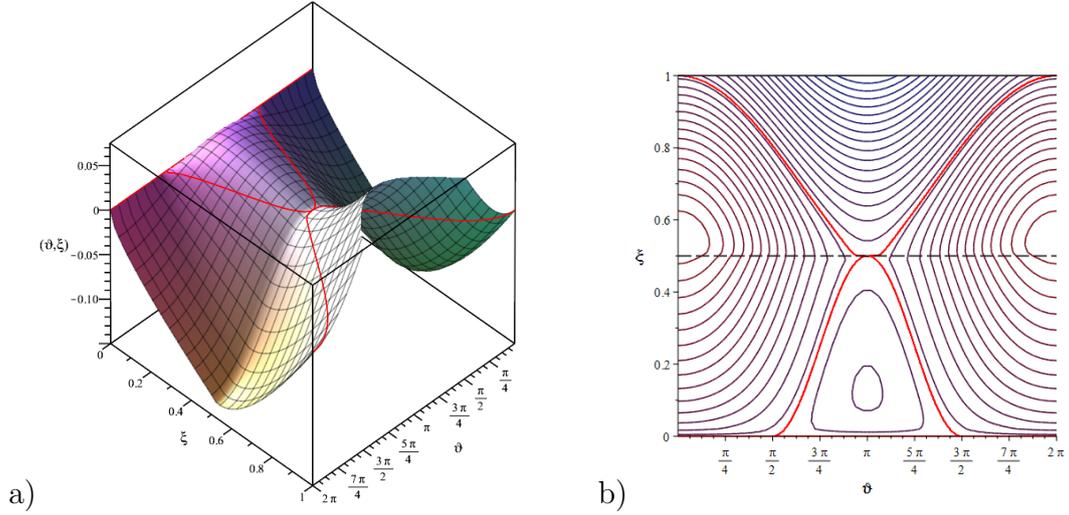

Figure 13: Coexistence of both bifurcation mechanisms for critical energy, detuning, and forcing amplitude of $\tilde{\xi} = 1$, $\sigma^* = f^* = 1.28$, respectively; a) 3D plot of $C(\nu, \xi)$, b) 2D projection of $C(\nu, \xi)$ on the $(\nu, \xi)$ plane.

The transition boundaries that correspond to bifurcations of type II for multiple $\tilde{\xi}$ values are shown in Fig. 14. As one can see, right and the left branches of the transition boundary correspond to the saddle and maximum mechanisms, respectively. The boundary associated with the maximum mechanism overlaps the left branch of the transition boundary for $\tilde{\xi} = 1/2$. In this case, the transition boundaries of both types overlap, as expected.

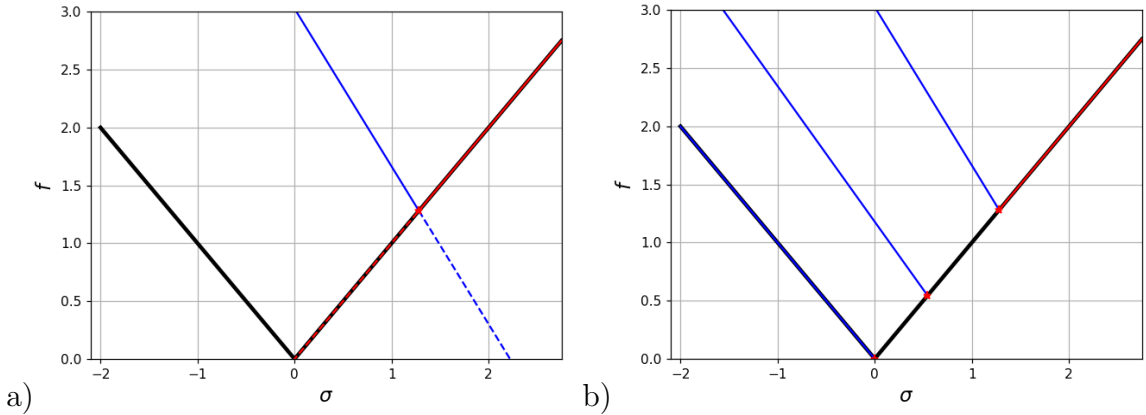

Figure 14: Transition boundaries that corresponds to bifurcation of type I (black line) and II. The right and left branches of the transition boundary associated with type II bifurcation correspond to the maximum (blue) and saddle (red) mechanisms, respectively. The intersection points of both mechanisms $(\sigma^*, f^*)$ are marked by a red star (Eq. (26)). Dashed line means that this escape mechanism is over-ruled by the other mechanism. a) for $\tilde{\xi} = 1$, b) from the bottom up: $\tilde{\xi} = 0.5, 0.75, 1$.

Bifurcation of type I corresponds to the transition boundary marked in black. Crossing this transition boundary towards the basin of HVI-regime involves a discrete increase in the response energy $\tilde{\xi}$. The energy level reached after crossing this transition boundary for a given detuning value $\tilde{\xi}^+(\sigma)$ is obtained by the inverse relation of the implicit expression shown in Eq. (27), where $\tilde{\xi}^+$ is the maximal response energy obtained immediately after crossing the transition boundary. Plus and minus signs in Eq. (27) correspond to the right and left branches (i.e. saddle and maximum mechanisms), respectively. For large negative detuning values, the energy increase tends to the following asymptotic value: $\lim_{\sigma \to \infty} \tilde{\xi}^+(\sigma) \approx 0.5543$.

$$\sigma \pm \frac{2}{\epsilon a_1(\tilde{\xi}^+)} \left( (1 + \epsilon\sigma) J(\tilde{\xi}^+) - \tilde{\xi}^+ \right) = 0 \qquad (27)$$

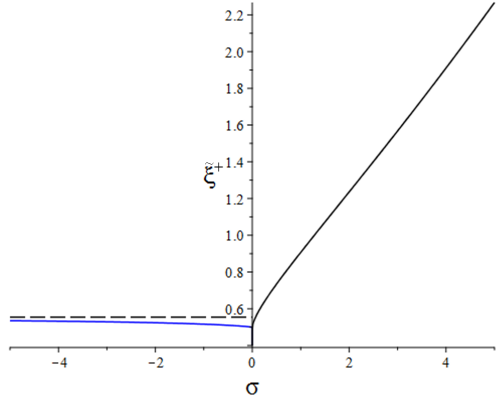

Figure 15: The maximal energy level $\tilde{\xi}^+(\sigma)$ obtained for a given detuning after crossing the transition boundary that corresponds to type II bifurcation for $\tilde{\xi} = 1$ towards the basin of the HVI-regime. Solid black and blue lines correspond to the right and left branches of the transition boundary, that correspond to bifurcation through the saddle and maximum mechanisms, respectively. Dashed black line corresponds to the limiting value for large negative detuning: $\lim_{\sigma \to \infty} \tilde{\xi}^+(\sigma) \approx 0.5543$.

As one can see in Fig. 15, crossing through the right branch of the transition boundary (bifurcation through the saddle mechanism) involves a much more drastic increase in the response energy $\tilde{\xi}$ in comparison to the left branch (bifurcation through the maximum mechanism). This effect can be also seen in the frequency response curve of the system, that corresponds to the explicit expression in Eq. (28) and plotted in Fig. 15. Since the response displacement in not informative due to the existence of rigid barriers at $|q| = 1$, the frequency response is given in terms of the energy of the response.

$$f - f_m(\sigma|\tilde{\xi}) = 0 \qquad (28)$$

Here, $f_m(\sigma|\tilde{\xi})$ is taken from Eq. (22). As one can see in Fig. 16, the energy jumps associated with crossing the transition boundary associated with type I bifurcations, and thus take place for $f = |\sigma|$. The increase measures are correlated with the values obtained by Eq. (27) and shown in Fig. 15.

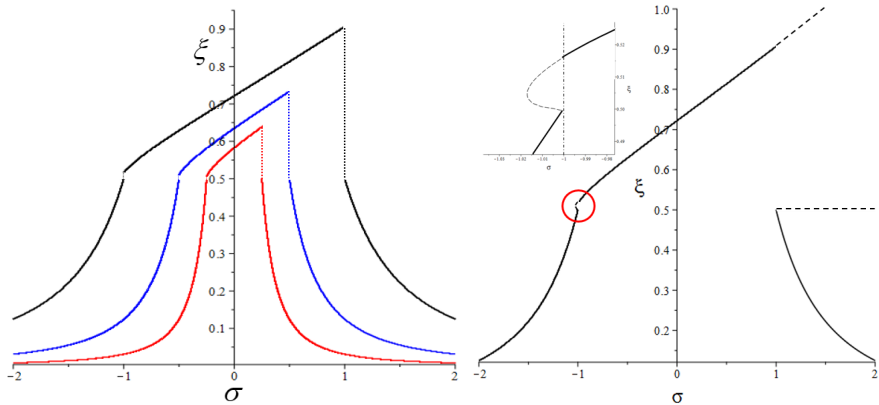

Figure 16: a) Frequency response curves for $f = 0.25, 0.5, 1$, colored in red, blue, and black, respectively. Dashed lines correspond to non-smooth energy jumps associated with crossing the transition boundary that correspond to bifurcation of type II. Left and right jumps correspond to crossing the left and right branches of the transition boundary, i.e. bifurcation through the maximum mechanism and saddle mechanism, respectively. b) Frequency response for $f = 1$ and zoom-in around the energy jump associated with crossing the left branch of the transition boundary. Dashed lines correspond to unstable branches of the frequency response.

Finally, Eq. (22)-(25) are used to plot function $\tilde{\xi}(\sigma, f)$ that maps the forcing parameters to the resulting maximal response energy reached by the HVI-oscillator. Contour lines $\tilde{\xi}(\sigma, f) = const$ represent iso-energy contours above the forcing parameters plane $\sigma - f$. The representation shown in Fig. 17 gives a full perspective on the response regimes, response energy levels, and energy absorption capabilities of the oscillator over the entire space of monochromatic harmonic excitations.

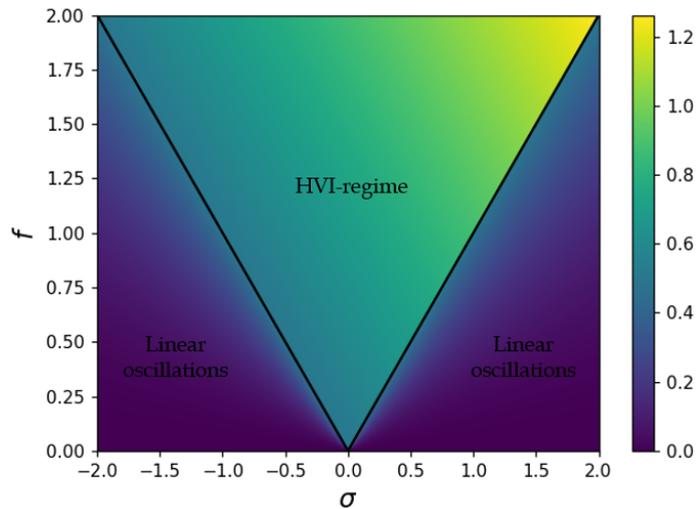

Figure 17: Plot of the maximal response energy $\tilde{\xi}(\sigma, f)$ over the forcing parameters plane. Black solid lines correspond to the transition boundary associated with type I bifurcation.

## 5. Numerical verification

Further exploration verifies the analytical results using numerical integration of the system's equation of motion, Eq. (1). First, the bifurcation of type I is investigated both in time domain and in the forcing parameters space. Both bifurcation mechanisms are illustrated in Fig. 18 and Fig. 19[1], respectively. As one can see in Fig. 18, minor increase in the forcing amplitude leads to a small increase in the maximal transient in the system's response. On the other hand, in Fig. 19, the same change measure in the forcing amplitude results in a dramatic increase in the response energy. Then, type II bifurcation for energy level $\tilde{\xi} = 1$ (chosen arbitrarily) through both bifurcation mechanisms is illustrated in Fig. 20-Fig. 21, and demonstrates similar behaviour. Finally, the theoretical boundaries associated with transition between linear oscillations and HVI-regime ($\tilde{\xi} = 1/2$), and transition to HVI-regime with maximal response energy of $\tilde{\xi} = 1$ are verified numerically in Fig. 22.

---

[1]As mentioned above, for type I bifurcation energy values above the bifurcation energy level, i.e. $\xi > \tilde{\xi}$ are irrelevant. Thus, any bifurcation of type II for any critical energy value $\tilde{\xi}$ is defined as type I bifurcation through maximum mechanism.

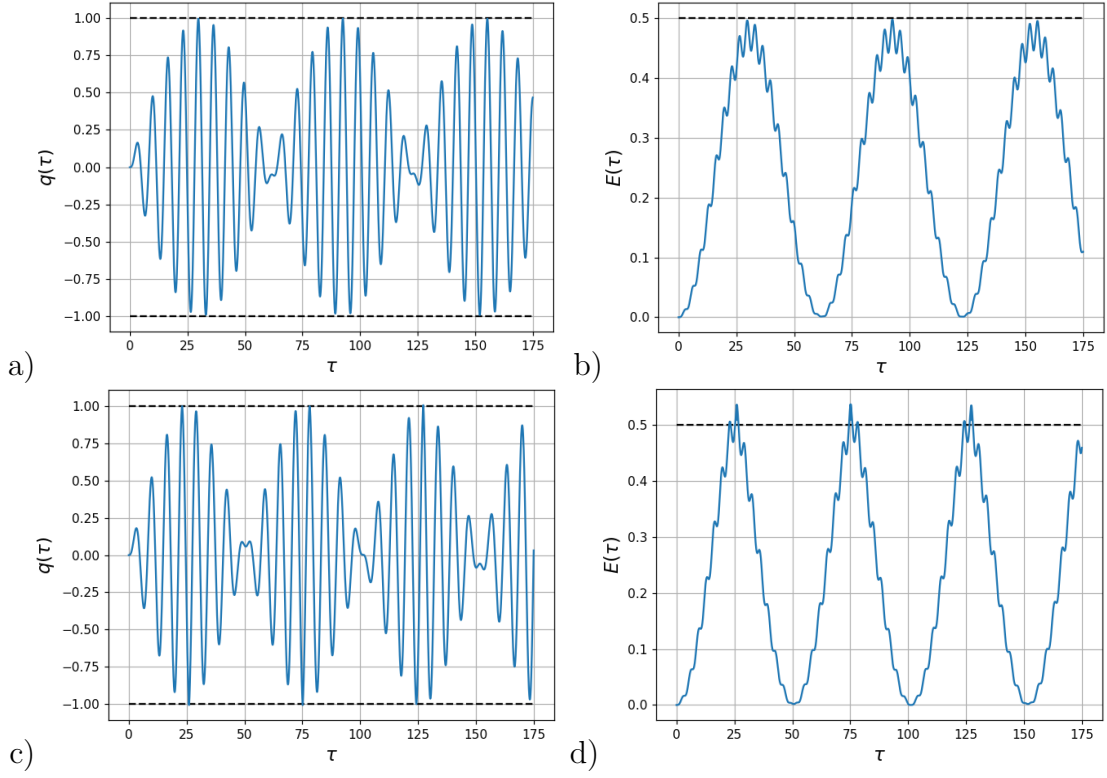

Figure 18: Bifurcation of type I ($\tilde{\xi} = 1/2$) through the maximum mechanism, for $\sigma = -1$ and a-b) $f = 0.99$, c-d) $f = 1.0$. a,c) Displacement responses, b,d) energy responses. Dashed black lines correspond to the energy bifurcation value, i.e. $\tilde{\xi} = 1/2$.

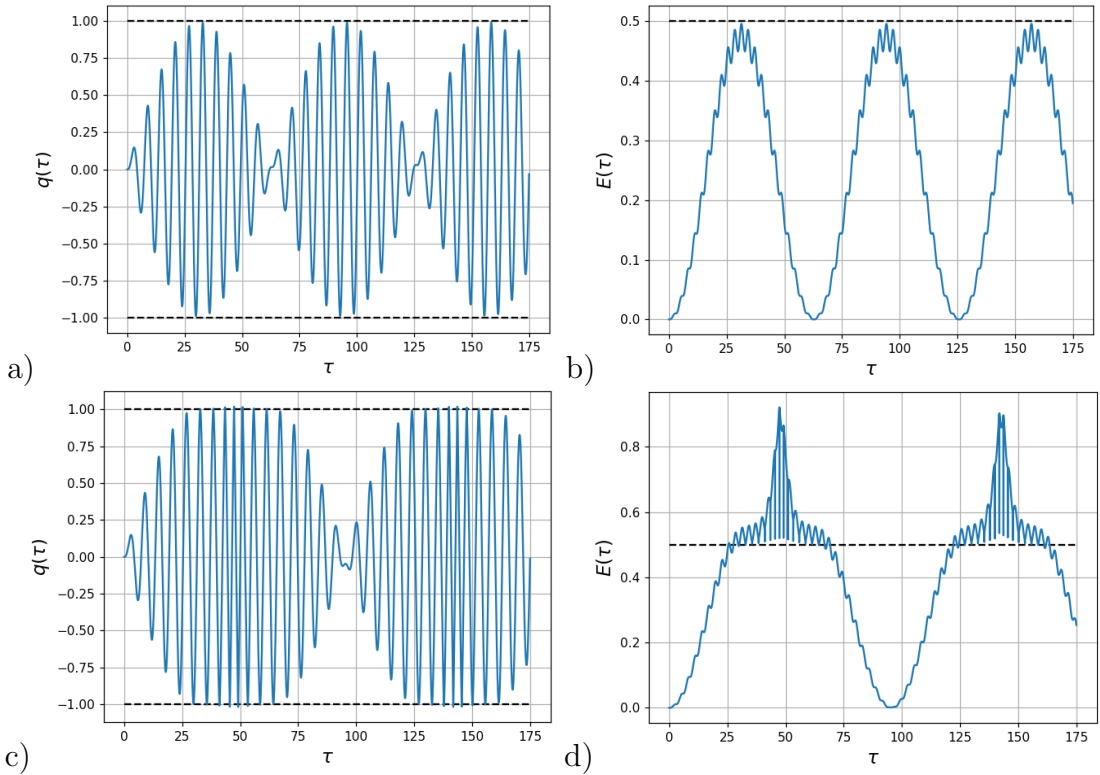

Figure 19: Bifurcation of type I ($\tilde{\xi} = 1/2$) through the saddle mechanism[1], for $\sigma = -1$, a-b) $f = 0.99$, c-d) $f = 1.0$. a,c) Displacement responses, b,d) energy responses. Dashed black lines correspond to the bifurcation energy value, i.e. $\tilde{\xi} = 1/2$.

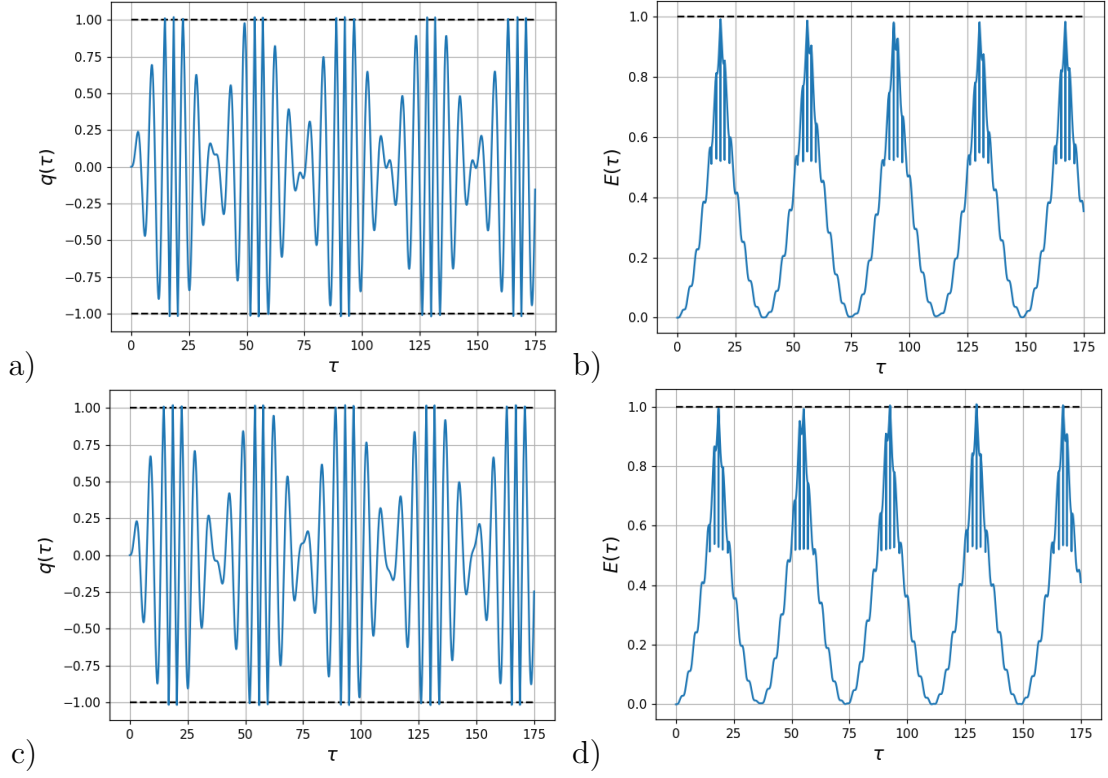

Figure 20: Bifurcation of type II ($\tilde{\xi} = 1/2$) through the maximum mechanism, for $\sigma = 1$ i.e. $f_m(\sigma = 1, \tilde{\xi} = 1) = 1.6637$, a-b) $f = 1.6$, c-d) $f = 1.7$. a,c) Displacement responses, b,d) energy responses. Dashed black lines correspond to the bifurcation energy value, i.e. $\tilde{\xi} = 1$.

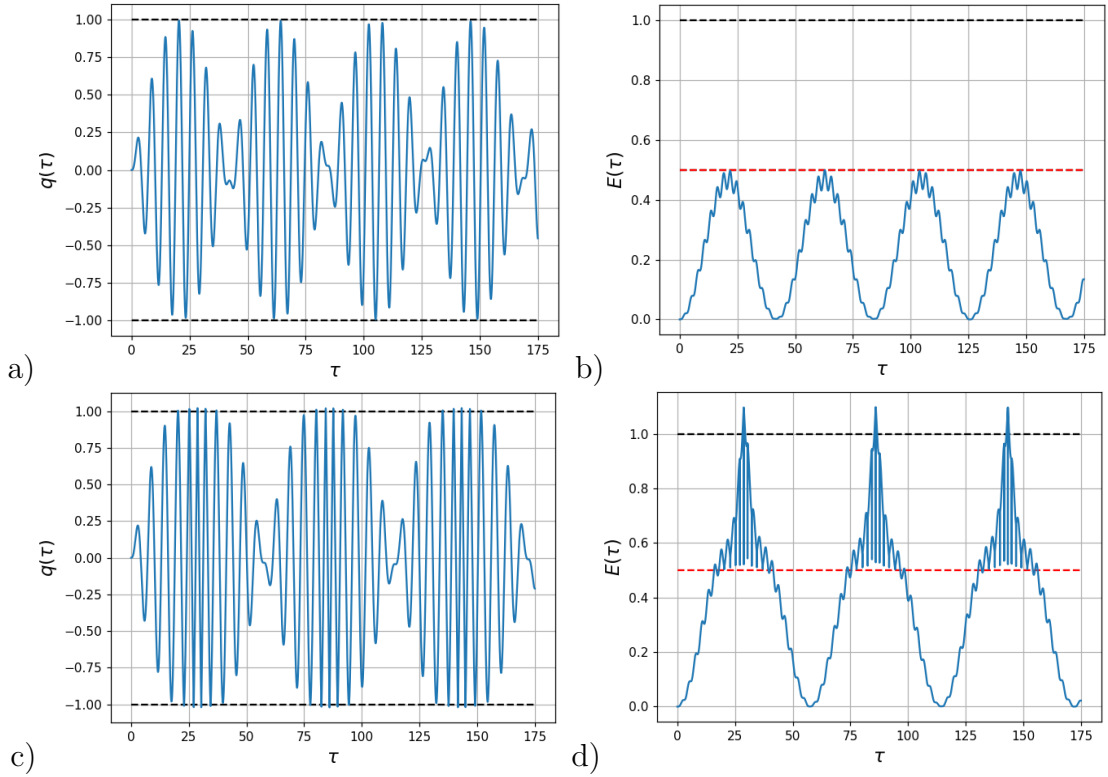

Figure 21: Bifurcation of type II for critical energy level of $\tilde{\xi} = 1$ through the saddle mechanism, for $\sigma = 1.5$, i.e. $f_s(\sigma = 1, \tilde{\xi} = 1) = 1.5$; a-b) $f = 1.49$, c-d) $f = 1.51$. a,c) Displacement responses, b,d) energy responses. Dashed black and red lines correspond to the bifurcation energy value, i.e. $\tilde{\xi} = 1$, and the saddle point at $\xi = 1/2$, respectively.

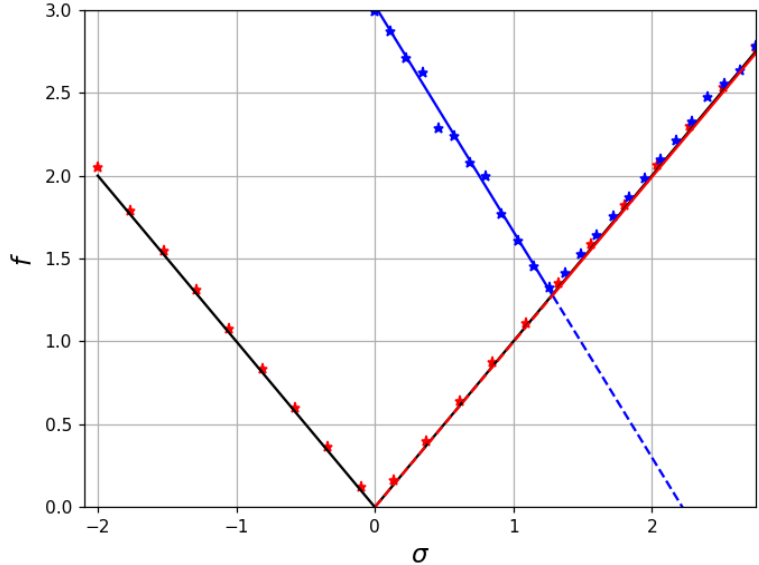

Figure 22: The transition boundaries that correspond to bifurcations of type I (black line) and II for critical energy of $\tilde{\xi} = 1$ (solid blue and red lines) and numerical verifications. The blue and red solid lines correspond to type II bifurcation through the maximum and saddle mechanisms, respectively. Dashed line means that the bifurcation mechanism is over-ruled by the other. Numerical verifications for the transition boundaries of bifurcations of type I and II are marked by red and blue stars, respectively. Numerical verification was obtained by integrating the system's equation of motion (Eq. (1)).

## 6. Concluding remarks

The main result of this work is finding analytical predictions for the transient maximal energy response reached in the system of a harmonically-forced hybrid vibro-impact oscillator for a given set of forcing parameters. Since the latter is a strongly non-linear system, traditional perturbation-based methods are inapplicable for the aforementioned task. Hence, canonical formalism was used to describe the dynamics in terms of action-angle variables. The generalized basis functions of the HVI-oscillator were derived. The slow flow dynamics in the vicinity of primary resonance was described by a cylindrical reduced resonance manifold, and two types of bifurcations were identified: the former is associated with the transition from linear oscillations to the hybrid vibro-impact regime and vice versa, and the latter with reaching a chosen transient energy level. Both bifurcations types correspond to transition boundaries in the forcing parameters plane, and they overlap when the critical energy value equals the energy levels that correspond to the birth of the HVI-regime, i.e. $\tilde{\xi} = 1/2$. The sharp minimum in the transition boundary curves is associated with the competition between two topological mechanisms of the LPT: the maximum mechanism and the saddle mechanism. These mechanisms are associated with the left and right branches of the transition boundary associated with type II bifurcation, respectively. This pattern is universal for escape problems under periodic forcing. Frequency response curves in terms of the transient energy of the system's forced response were obtained. Finally, an analytical prediction of the resulting maximal transient energy level over the forcing parameters space was obtained analytically. All analytical predictions were verified numerically. The results also point out the advantage of the energy-based methods for the prediction of response regimes in strongly non-linear dynamical systems, and action-angle variables in particular.

The aforementioned results give a full perspective about the vibration mitigation performances of an HVI-based passive energy absorber (HVI-NES [9]) and can be used for accurate optimization of the absorber design parameters, prediction of its absorption rate for a given set of excitation properties, and hence, the resistance of the primary structure. We learn that the saddle mechanism can lead to catastrophic consequences in HVI-like engineering systems, such as fluid sloshing in liquid storage tanks, flow-induced vibration, and impacts between nuclear fuel rods and grids under axial coolant flow, that might potentially result in a drastic increase in the response energy due to a slight uncontrolled change in the forcing parameters, and finally lead to mechanical wear and failure. On the other hand, for HVI-NES the saddle mechanism allows significant enhancement in the absorption rate and thus improves vibration protection capabilities. Those findings and insights are to be considered in the design and optimization process of such a vibration absorber.

**Acknowledgements** The author is grateful to Prof. V. Pilipchuk for useful discussions.

**Funding** M. Farid has been supported by the Fulbright Program, the Israel Scholarship Education Foundation (ISEF), Jean De Gunzburg International Fellowship, the Israel Academy of Sciences and Humanities, the Yitzhak Shamir Postdoctoral Scholarship of the Israeli Ministry of Science and Technology, the PMRI – Peter Munk Research Institute - Technion, and the Israel Science Foundation Grant 1696/17.

## Appendix A. Detailed derivation of action-angle transformation

In this appendix we bring the detailed derivation of expressions used in the canonical transformation of the dynamical system presented. The relation between displacement and momentum of the HVI-oscillator is given by the following standard conservation rule:

$$\frac{1}{2}\left(q^2 + \dot{q}^2\right) = E \to \dot{q} = \pm\sqrt{2E - q^2} \tag{A.1}$$

According to Eq. (8), the action variable is calculated as follows:

$$I(E) = \frac{2}{\pi} \int_0^1 \sqrt{2E - q^2}\,dq = \frac{1}{\pi}\left(\sqrt{2E-1} + 2E \arctan\left(\frac{1}{\sqrt{2E-1}}\right)\right) \tag{A.2}$$

The frequency is calculated as follows:

$$\frac{1}{\omega(E)} = \frac{\partial I}{\partial E} = \frac{2}{\pi}\arctan\left(\frac{1}{\sqrt{2E-1}}\right) \to \omega(E) = \frac{\pi}{2}\left(\arctan\left(\frac{1}{\sqrt{2E-1}}\right)\right)^{-1} \tag{A.3}$$

For the well-known case of linear oscillator, we know that the frequency equals to unity and the action equals to the energy of the oscillator for all energy values, i.e. $\omega(E) = 1, I(E) = E$. Hence, to merge both cases, the second expression in Eq. (14) was introduced. The angle variable is derived according to Eq. (8):

$$\begin{aligned}\theta &= \frac{\partial}{\partial I}\int_0^q \sqrt{2(E-U(x))}\,dx = \omega(E)\frac{\partial}{\partial E}\int_0^q \sqrt{2(E-U(x))}\,dx \\ &= \frac{\omega(E)}{\sqrt{2E}}\int_0^q \frac{1}{\sqrt{1-\frac{U(x)}{E}}} = \omega(E)\arcsin\left(\frac{x}{\sqrt{2E}}\right)\bigg|_0^q = \omega(E)\arcsin\left(\frac{x}{\sqrt{2E}}\right)\end{aligned} \tag{A.4}$$

Inverting Eq. (A.4) yields the following solution in terms of energy and angle:

$$q(E, \theta) = \sqrt{2E}\sin\left(\frac{\theta}{\omega(E)}\right) \tag{A.5}$$

In accordance to Eq.(13), the coefficient of the first term in the Fourier series of Eq. (16) is required to obtain the RM of the system:

$$q_1(E) = -\frac{i}{2}a_1(E), \quad a_1(E) = \frac{1}{\pi}\int_{-\pi}^{\pi}\sqrt{2E}\sin\left(\frac{\theta}{\omega}\right)\sin\theta\,d\theta = \frac{\sqrt{2E}2\omega^2\sin\left(\frac{\pi}{\omega}\right)}{\pi(\omega^2 - 1)} \tag{A.6}$$

# List of Abbreviations

| | |
|---|---|
| AA | action-angle (variables) |
| HVI | Hybrid vibro-impact (regime, oscillator, NES) |
| LPT | Limiting phase trajectory |
| NES | nonlinear energy sink |
| PEA | Passive energy absorber |
| RM | Resonance manifold |
| TMD | Tuned mass damper |
| TRC | Transient resonance capture |

# List of Symbols

| | |
|---|---|
| $2d$ | The channel's length |
| $\bar{\tau}(\tau), \bar{e}(\tau)$ | Basis functions of the VI-oscillator |
| $\bar{F}(t)$ | Dimensional time-dependant forcing |
| $\delta$ | The Dirac delta function |
| $\epsilon$ | Small parameter |
| $\kappa$ | Restitution coefficient |
| $\nu$ | Phase variable |
| $\omega$ | Frequency of oscillations of the oscillator |
| $\sigma$ | Detuning parameter |
| $\sigma^*, f^*$ | Detuning and critical forcing amplitude associated with coexistence of both 'saddle mechanism' and 'maximum mechanism' |
| $\tilde{\xi}$ | The energy level |
| $\tilde{\xi}^+$ | The energy level associated with crossing the transition boundary |
| $\xi_s$ | The transient energy level associates with the saddle point of the RM |
| $a_k$ | The coefficient of the $k^{th}$ term in Fourier series of solution $q(I, \theta)$ |
| $C(\nu, \xi)$ | The expression describes the resonance manifold |
| $E, \xi$ | Instantaneous and averaged energy of the oscillator |
| $f$ | Scaled forcing amplitude |
| $F, \Omega$ | Non-dimensional forcing amplitude and frequency |
| $f_m, f_s$ | The critical forcing amplitudes associated with the occurrence of bifurcation through through maximum and saddle mechanisms, respectively |
| $g(\tau), g'(\tau), \beta$ | Basis functions of the HVI-oscillator and energy-related parameter |
| $H$ | Integral of motion/conservation law |
| $H_0$ | Initial conditions-related value of the integral of motion, H |
| $i$ | Unit imaginary number |
| $I, \theta$ | Action and angle variables |
| $J$ | Averaged action variable |
| $k$ | Dimensional linear stiffness coefficient |
| $p$ | Momentum of the oscillator |
| $t, \tau$ | Dimensional and non-dimensional time scales |
| $U(q)$ | Equivalent potential energy function |
| $x, q$ | Dimensional and non-dimensional displacement of the oscillator |